\documentclass[10pt,reqno]{amsart}

\usepackage{amssymb}
\usepackage{amsthm}
\usepackage{amsmath}

\newtheorem{thm}{Theorem}

\newtheorem{prop}[thm]{Proposition}

\newtheorem{knownthm}{Theorem}

\newtheorem{knownthmp}{Theorem}

\theoremstyle{definition}

\newcommand{\closure}{\overline}

\newcommand{\round}{\partial}
\newcommand{\CC}{\widehat{\mathbb{C}}}
\newcommand{\C}{\mathbb{C}}
\newcommand{\D}{\mathbb{D}}

\newcommand{\R}{\mathbb{R}}
\newcommand{\B}{\mathcal{B}}
\newcommand{\A}{\mathcal{A}}
\renewcommand{\S}{\mathcal{S}}
\newcommand{\no}{\noindent}
\newcommand{\dstyle}{\displaystyle}
\renewcommand{\Re}{\textup{Re}\,}
\renewcommand{\Im}{\textup{Im}\,}

\newcommand{\arctanh}{\textup{arctanh}\,}

\title[]{L\"owner chains with complex leading coefficient}
\author[]{IKKEI HOTTA}
\address{Division of Mathematics, Graduate School of Information Sciences, Tohoku University, 6-3-09 Aramaki-Aza-Aoba, Aoba-ku, Sendai 980-8579, Japan, phone : +81-22-795-4636}
\email{ikkeihotta@ims.is.tohoku.ac.jp}
\urladdr{http://www.ims.is.tohoku.ac.jp/{\textasciitilde}ikkeihotta}
\subjclass[2000]{Primary 30C55, 30C80, Secondary 30C45,30C62}
\keywords{L\"owner(Loewner) chain, quasiconformal mapping, univalent function.}
\date{Oct. 26, 2009}

\begin{document}

\begin{abstract}
In this paper we confirm that several crucial theorems due to Pommerenke and Becker for the theory of L\"owner chains work well without normalization on the complex-valued first coefficient.
As applications of those considerations, some new univalent and quasiconformal extension criteria are given in the last section.
\end{abstract}

\maketitle


\section{\bf Introduction}

Let $\C$ and $\CC$ denote the complex plane and the Riemann sphere respectively, and let $\D_{r} = \{z \in \C : |z| < r \leq 1\}$ and $\D = \D_{1}$.
We denote by $\A$ the class of functions $f(z)$ normalized so that $f(0)=f'(0)-1=0$ which are analytic on $\D$  and $\S$ the subclass of $\A$ whose members are univalent on $\D$.

Let $f_{t}(z) = f(z,t) = \sum_{n=1}^{\infty}a_{n}(t)z^{n},\,a_{1}(t) \neq 0,$ be a function defined on $\D \times [0,\infty)$ and analytic in $\D$ for each $t \in [0,\infty)$, where $a_{1}(t)$ is a complex-valued, locally absolutely continuous function on $[0,\infty)$ and $\lim_{t \to \infty}|a_{1}(t)| = \infty$.
Then $f_{t}(z)$ is called a \textit{L\"owner chain} if $f_{t}(z)$ satisfies the following conditions;

\def\labelenumi{\arabic{enumi}.}
\begin{enumerate}
\item $f_{t}(z)$ is univalent in $\D$ for each $t \in [0,\infty)$,
\item $f_{s}(\D) \subsetneq f_{t}(\D)$ for $0 \leq s < t < \infty$.
\end{enumerate}
If $a_{1}(t) = e^{t}$, then we say that $f_{t}(z)$ is a \textit{standard L\"owner chain}.

It is known that if $f_{t}(z)$ is a L\"owner chain then $f_{t_{n}}(\D) \to f_{t_{0}}(\D)$ if $t_{n} \to t_{0} \in [0,\infty)$ and $f_{t_{n}}(\D) \to \C$ if $t_{n} \to \infty$ in the sense of kernel convergence with respect to the origin.
However, the converse is not true in general (see for instance \cite[pp.136--138]{Conway:1995} and \cite[pp.94--97]{GrahamKohr:2003}).

L\"owner chains and several related theorems serve as a powerful tool for the theory of univalent functions.
In those researches, it seems that many authors deal with the case when the first coefficient $a_{1}(t)$ of a L\"owner chain is a real-valued function.
Mostly, they only require the condition $a_{1}(t) \neq 0$ and it is not clear that $a_{1}(t)$ is either real- or complex-valued.
Certainly, the fact is known that $a_{1}(t)$ can be taken to be complex-valued (e.g. \cite{Becker:1976}) and actually it is used by some authors as a tool (e.g. \cite{Pom:1965}, \cite{Ruscheweyh:1976}, \cite{Hotta:b}).
However, usually there is no explicit mention that the first coefficient of L\"owner chains is admitted to be a complex-valued function.

One of our purposes is to give a precise statement and a proof to our key theorems for L\"owner chains with complex first coefficient, and to confirm that those theorems work well without normalization.
In consequence, the obscurity which is mentioned above will be completely removed.

A L\"owner chain without normalization on the first derivative has some new applications to univalence and quasiconformal extension criteria. 
Another purpose of this paper is to derive those new criteria as a benefit of our considerations.
We discuss it in the last section.


\

\section{\bf Related theorems for L\"owner chains}

The following necessary and sufficient condition for a standard L\"owner chain is known;

\def\labelenumi{\roman{enumi})}

\begin{knownthm}[\cite{Pom:1965},\cite{Pom:1975}]\label{pom}
Let $0 < r_{0} \leq 1$.
Let $h(z,t) = e^{t}z + \sum_{n=2}^{\infty}c_{n}(t)z^{n}$ be a function defined on $\D \times [0,\infty)$ and analytic in $\D$ for each $t \in [0,\infty)$.
Then the function $h(z,t)$ is a standard L\"owner chain if and only if the following two conditions are satisfied;
\begin{enumerate}
\item The function $h(z,t)$ is analytic in $z \in \D_{r_{0}}$ for each $t \in [0,\infty)$, absolutely continuous in $t \in [0,\infty)$ for each $z \in \D_{r_{0}}$ and satisfies
\begin{equation}\label{pomineq}
|h(z,t)| \leq K_{0} e^{t} \hspace{20pt} (z \in \D_{r_{0}},\,t \in [0,\infty))
\end{equation}
for some positive constants $K_{0}$.
\item There exists a function $p(z,t)$ analytic in $z \in \D$ for each $t \in [0,\infty)$ and measurable in $t \in [0,\infty)$ for each $z \in \D$ satisfying
$$
\Re p(z,t) > 0 \hspace{20pt} (z \in \D,\,t \in [0,\infty))
$$
such that
\begin{equation} \label{LDE}
\dot{h}(z,t) =z h'(z,t) p(z,t)  \hspace{20pt} (z \in \D_{r_{0}},\,\textup{a.e.}~t \in [0,\infty))
\end{equation}
where $\dot{h} = \round h /\round t$ and $h' = \round h/\round z$.
\end{enumerate}
\end{knownthm}

We will show first that Theorem \ref{pom} can be generalized for a L\"owner chain which has the complex-valued first coefficient as the following form;

\def\labelenumi{\roman{enumi}')}

\begin{knownthmp}\label{pom2}
Let $0 < r_{1} \leq 1$.
Let $f(z,t) = \sum_{n=1}^{\infty}a_{n}(t)z^{n},\,a_{1}(t) \neq 0,$ be a function defined on $\D \times [0,\infty)$ and analytic in $\D$ for each $t \in [0,\infty)$, where $a_{1}(t)$ is a complex-valued, locally absolutely continuous function on $[0,\infty)$ and $\lim_{t \to \infty}|a_{1}(t)| = \infty$.
Then the function $f(z,t)$ is a L\"owner chain if and only if the following conditions are satisfied;
\begin{enumerate}
  \item The function $f(z,t)$ is analytic in $\D_{r_{1}}$ for each $t \in [0,\infty)$, locally absolutely continuous in $[0,\infty)$ for each $z \in \D_{r_{1}}$ and
\begin{equation}\label{pom2ineq}
  |f(z,t)| \leq K_{1} |a_{1}(t)| \hspace{20pt} (z \in \D_{r_{1}},\,\textup{a.e.}~t \in [0,\infty))
\end{equation}
for some positive constants $K_{1}$.
  \item There exists a function $p(z,t)$ analytic in $\D$ for each $t \in [0,\infty)$ and measurable in $[0,\infty)$ for each $z \in \D$ satisfying
  $$
\Re p(z,t) > 0 \hspace{20pt} (z \in \D,\,t \in [0,\infty))
$$
such that
\begin{equation}\label{pom2LDE}
\dot{f}(z,t) =z f'(z,t) p(z,t)  \hspace{20pt} (z \in \D_{r_{1}},\,\textup{a.e.}~t \in [0,\infty))
\end{equation}
where $\dot{f} = \round f /\round t$ and $f' = \round f/\round z$.
\end{enumerate}
\end{knownthmp}

\begin{proof}
Let $f(z,t) = \sum_{n=1}^{\infty}a_{n}(t)z^{n}$ be a L\"owner chain, where $a_{1}(t) \neq 0$ is a complex-valued locally absolutely continuous function on $[0,\infty)$.

First we set $\lambda(t) := -\arg a_{1}(t)$ which is locally absolutely continuous with respect to $t$
and define
\begin{equation}\label{g(z,t)}
g(z,t) = \sum_{n=1}^{\infty}b_{n}(t)z^{n} := f(e^{i\lambda(t)}z,t).
\end{equation}
This yields
$$
  \frac{\dot{g}(z,t)}{zg'(z,t)} = \frac{\dot{f}(e^{i\lambda(t)}z,t)}{e^{i\lambda(t)}zf'(e^{i\lambda(t)}z,t)} + i \lambda'(t).
$$
It follows that $\Re\dot{g}(z,t)/z g'(z,t) > 0 $ if and only if $\Re\dot{f}(z,t)/zf'(z,t)>0$ for all $z \in \D$ and almost all $t \in [0,\infty)$.
We remark that $b_{1}(t) = |a_{1}(t)|$ is a strictly increasing positive function in $t \in [0,\infty)$, and the inverse function $b_{1}^{-1}$ is defined in $[|a_{1}(0)|, \infty)$ and maps this interval onto $[0,\infty)$.

Next, let
\begin{equation}\label{h(z,t)}
h(z,t) 
:= 
\frac{1}{|a_{1}(0)|}
g(z,b_{1}^{-1}(|a_{1}(0)|e^{t})).
\end{equation}
It follows immediately from these reparametrizations that $f(z,t)$ is a L\"owner chain if and only if $h(z,t)$ is a standard L\"owner chain.
We note that therefore $h$, namely $b_{1}^{-1}$ also, is absolutely continuous with respect to $t$ by Theorem \ref{pom}.

We also have
$$
\frac{\dot{h}(z,t)} {zh'(z,t)} 
=
\frac{\dot{g}(z,b_{1}^{-1}(|a_{1}(0)|e^{t}))}{zg'(z,b_{1}^{-1}(|a_{1}(0)|e^{t}))}
\cdot (b_{1}^{-1}(|a_{1}(0)|e^{t}))'.
$$
from \eqref{h(z,t)}.
Since $(b_{1}^{-1}(t))' > 0$ for almost all $t \in [0,\infty)$, $\Re\dot{h}(z,t)/z h'(z,t) > 0 $ if and only if $\Re\dot{g}(z,t)/z g'(z,t) > 0 $ for all $z \in \D$ and $t \in [0,\infty)$. 
We can see easily that $h(z,t)$ satisfies the condition \eqref{pomineq} if and only if $f(z,t)$ satisfies \eqref{pom2ineq} with $K_{1} = K_{0}/|a_{1}(0)|^{2}$, and the other properties of the sufficient part of Theorem \ref{pom} are preserved by the reparametrizations \eqref{g(z,t)} and \eqref{h(z,t)} with $r_{0} = r_{1}$.
Consequently, all the necessary and sufficient conditions of Theorem \ref{pom2} follows from Theorem \ref{pom}.
\end{proof}

We remark that a similar argument about normalization of L\"owner chains as above can be found in \cite{Becker:1976}. 
In fact, $p(z,t)$ in \eqref{pom2LDE} is normalized by
$$
p^{*}(e^{i\beta(t)}z, \alpha(t)) = \frac{1}{\Re p(0,t)}\left[p(z,t) -i \Im p(0,t)\right]
$$
where
$$
\alpha(t) = \int_{0}^{t}\Re p(0,\tau) d\tau,
\hspace{20pt}
\beta(t) = \int_{0}^{t}\Im p(0,\tau) d\tau,
$$
which is equivalent to $f'(0,t) = \exp\{\alpha(t) + i \beta(t)\}$.

Theorem \ref{pom} appears in \cite[p.4]{Millermocanu:2000} without proof, and the similar change of variables as \eqref{g(z,t)} and \eqref{h(z,t)} is in \cite[p.95]{GrahamKohr:2003}. 
However, it is not clear that whether $a_{1}(t) = f'(0,t)$ can be taken to be complex-valued or not.

The following theorem which is essentially same as Theorem \ref{pom2} is often used to show univalence for an analytic function, apart from the theory of L\"owner chains;

\begin{knownthm}[\cite{Pom:1965}]
Let $0 < r_{0} \leq 1$.
Let $f(z,t) = a_{1}(t)z + \sum_{n=2}^{\infty}a_{n}(t)z^{n},\,a_{1}(t) \neq 0$,\, be analytic for each $t \in [0,\infty)$ in $\D_{r_{0}}$ and locally absolutely continuous in $[0,\infty)$, locally uniformly with respect to $\D_{r_{0}}$, where $a_{1}(t)$ is a complex-valued function on $[0,\infty)$.
For almost all $t \in [0,\infty)$ suppose
$$
\dot{f}(z,t) =z f'(z,t) p(z,t)  \hspace{20pt} (z \in \D_{r_{0}},\,t \in [0,\infty))
$$
where $p(z,t)$ is analytic in $\D$ and satisfies $\Re p(z,t)>0,\,z \in \D$.
If $|a_{1}(t)| \to \infty$ for $t \to \infty$ and if $\{f(z,t)/a_{1}(t)\}$ forms a normal family in $\D_{r_{0}}$, then for each $t \in [0,\infty)$ $f(z,t)$ can be continued analytically in $\D$ and gives a univalent function.
\end{knownthm}

\begin{proof}
By using the previous reparametrization argument, it is enough to think about the case when $f(z,t)$ is a standard L\"owner chain.
Following the lines of the proof of \cite[Folgerung 3]{Pom:1965} one can obtain our assertion.
\end{proof}

Furthermore, we shall look at the next theorem which is due to Becker.
Here, a sense-preserving homeomorphism $f$ of $G \subset \C$ is called \textit{$k$-quasiconformal} if $f_{z}$ and $f_{\bar{z}}$, the partial derivatives in $z$ and $\bar{z}$ in the distributional sense, are locally integrable on $G$ and satisfy $|f_{\bar{z}}| \leq k|f_{z}|$ almost everywhere in $G$, where $k \in [0,1)$.

\begin{knownthm}[\cite{Becker:1972}]\label{Beckerthm}
Suppose that $ h_{t}(z) = h(z,t)$ is a standard L\"owner chain for which $p(z,t)$ in \eqref{LDE} satisfies the condition 
\begin{equation*}
\begin{array}{lrll}
p(z,t) \in U(k) &\!:=\!&
\dstyle
\left\{
w \in \C : \left|\frac{w-1}{w+1}\right| \leq k
\right\}\\[15pt]
&\!=\!&
\dstyle
\left\{
w \in \C: \left|w-\frac{1+k^{2}}{1-k^{2}}\right| \leq \frac{2k}{1-k^{2}}
\right\}
\end{array}
\end{equation*}
for all $z \in \D$ and $t \in [0,\infty)$.
Then $h(z,t)$ admits a continuous extension to $\closure{\D}$ for each $t \geq 0$ and the map $\hat{h}$ defined by
\begin{equation}\label{becker}
\hat{h}(z) =
\left\{
\begin{array}{ll}
\dstyle h(z,0), &\textit{if}\hspace{10pt}   |z| < 1, \\[5pt]
\dstyle h(\frac{z}{|z|},\log |z|), &\textit{if}\hspace{10pt}   |z| \geq 1,
\end{array} 
\right.
\end{equation}
is a $k$-quasiconformal extension of $h_{0}$ to $\C$.
\end{knownthm}

For the proof, see e.g. \cite{Becker:1980}.
The above theorem is also generalized for a L\"owner chain $f(z,t)$ with the complex first coefficient;

\setcounter{knownthmp}{2}
\begin{knownthmp}[\cite{Becker:1976}]\label{Beckerthm2}
Suppose that $f_{t}(z) = f(z,t)$ is a L\"owner chain for which $p(z,t)$ in \eqref{LDE} satisfies the condition $p(z,t) \in U(k)$ for all $z \in \D$ and $t \in [0,\infty)$, where $a_{1}(t) = f'(0,t)$ is a complex-valued function on $[0,\infty)$.
Then $f_{t}(z)$ admits a continuous extension to $\closure{\D}$ for each $t \geq 0$ and the map $\hat{f}$ given in \eqref{becker} is a $k$-quasiconformal extension of $f_{0}$ to $\C$.
\end{knownthmp}

For the proof, see \cite{Becker:1976}.
Indeed, it is enough to repeat the original proof in \cite{Becker:1980} with Theorem \ref{pom2} for $f(z,t)$.


\

\section{\bf Applications}
In this section we shall show some univalence and quasiconformal extension criteria as a consequence of the considerations in the previous section.

\

\subsection{Convex combinations}

A function $f \in \A$ is called \textit{convex} if $f$ is univalent and $f(\D)$ is a convex domain.
It is well known that $f \in \A$ is convex if and only if $f$ satisfies $\Re \{1 + zf''(z)/f'(z)\} > 0$ for all $z \in \D$.
The author showed in \cite{Hotta:2009a} that if $f$ is convex then $\alpha f (z) + (1-\alpha)zf'(z)$ with $\alpha \in [0,1]$ is univalent in $\D$.
This result is extended as following with the aid of Theorem \ref{pom2};

\begin{thm}\label{ccombi}
Let $\alpha$ be a complex number with $|2 \alpha - 1| \leq 1$.
If $f \in \A$ is a convex function, then the function
\begin{equation}\label{ccombi1}
 \alpha f(z) + (1 - \alpha) zf'(z)
\end{equation}
is univalent in $\D$.
\end{thm}

\begin{proof}
Let
$$
 f_{t}(z) = \alpha f(z) + e^{t} (1-\alpha) zf'(z).
$$
Then $\lim_{t \to \infty}|f_{t}'(0)| = \lim_{t \to \infty}|\alpha + (1 - \alpha)e^{t}| = \infty$.
Furthermore
$$
\frac{1}{p(z,t)} 
=
\frac{zf_{t}'(z)} {\dot{f_{t}}(z)}
= 
\frac{1}{e^{t}}\left(\frac{\alpha}{1-\alpha}\right) + 1 + \frac{zf''(z)}{f'(z)}.
$$
Since $\Re \alpha/(1-\alpha) >0$, it turns out that $f_{t}(z)$ is a L\"owner chain by the assumption and Theorem \ref{pom2}.
In particular $f_{0}$ is univalent in $\D$ which is our assertion.
\end{proof}

It can be shown that the function \eqref{ccombi1} is close-to-convex by using convolution technique (Li-Mei Wang, personal communications).
On the other hand, by using L\"owner's method we can give a simple proof for univalency of \eqref{ccombi1}.

\

\subsection{Spirallike functions}

A function $f \in \A$ is called \textit{$\alpha$-spirallike} and known to be univalent if $f$ satisfies
$$
\Re\left\{e^{i\alpha}\frac{zf'(z)}{f(z)} \right\}>0
$$
for a real number $\alpha \in (-\pi/2,\pi/2)$ in $\D$.
If $\alpha = 0$ then we say that $f$ is \textit{starlike}.
It is known \cite{Pom:1975} that the standard L\"owner chain
\begin{equation}\label{spiral1}
h_{t}(z) = e^{(1-i a)t} f(e^{i a t}z)
\end{equation}
with $a = \tan \alpha$ corresponds to an $\alpha$-spirallike function because it follows from calculations that
\begin{equation}\label{spiral2}
p(z,t)
=
\frac{\dot{h_{t}}(z)}{zh_{t}'(z)} 
=
ia + \frac{1}{\cos \alpha} 
\left(
e^{-i\alpha}\frac{f(e^{iat}z)} {e^{iat}zf'(e^{iat}z)}
\right)
\end{equation}
and therefore $\Re p(z,t)>0$ implies $\alpha$-spirallikeness of $f$. 
If we apply Theorem \ref{pom} and Theorem \ref{Beckerthm} to the standard L\"owner chain \eqref{spiral1}, then we obtain the next proposition.
Here, let us denote that $U(\alpha,k)$ the hyperbolic disk in the tilted half plane $\{z \in \C : \Re e^{i\alpha} z > 0\}$ centered at 1 with radius $\arctanh k,\,0 \leq k < 1$, i.e.,
$$
U(\alpha, k) =
\left\{
 w \in \C :
 \left|
  w - 
  \frac{1+e^{-2i\alpha}k^{2}}{1-k^{2}}
 \right|
 \leq \frac{2k\cos \alpha}{1-k^{2}}
\right\}.
$$
It is clear that $U(0,k) = U(k)$.

\begin{prop}
Let $\alpha \in (-\pi/2,\pi/2)$ and $k \in [0,1)$.
For $f \in \A$, if
$$
\frac{zf'(z)}{f(z)} \in U(\alpha,k)
$$
for all $z \in \D$, then $f$ has a k-quasiconformal extension to $\C$.
\end{prop}

The case when $\alpha = 0$ appears e.g. in \cite{Brown:1984}.

\begin{proof} 
Let $w = zf'(z)/f(z)$.
Then applying Theorem \ref{pom} and Theorem \ref{Beckerthm} to \eqref{spiral2}, it can be deduced that if $w$ satisfies
$$
 \left|
\frac1w -
  \frac{1+e^{2i\alpha}k^{2}}{1-k^{2}}
 \right|
 \leq \frac{2k\cos \alpha}{1-k^{2}}
$$
then, $f$ has a $k$-quasiconformal extension to $\C$.
The inequality implies $1/w \in U(-\alpha,k)$ which is equivalent to $w \in U(\alpha,k)$
\end{proof}

On the other hand, by constructing a L\"owner chain without normalization on the first derivative and applying Theorem \ref{pom2} and Theorem \ref{Beckerthm2} we have the following another result;

\begin{thm}\label{spiralthm}
Let $\alpha \in (-\pi/2,\pi/2)$ and $k \in [\,|\tan(\alpha/2)|,1)$.
For $f \in \A$, if
$$
e^{i \alpha} \frac{zf'(z)}{f(z)} \in U(k)
$$
for all $z \in \D$, then $f$ has a k-quasiconformal extension to $\C$.
\end{thm}

\begin{proof}
Let $c$ be a complex constant with $\Re c >0$.
If we set
$$
f_{t}(z) = e^{ct}f(z),
$$
then $\lim_{t \to \infty} |f_{t}'(0)| = \lim_{t \to \infty} |e^{ct}| = \infty$ since $\Re c>0$.
A calculation shows that
$$
\frac{1}{p(z,t)}
=
\frac{1}{c} \frac{zf'(z)}{f(z)}.
$$
Therefore we obtain our theorem if we put $c = e^{-i\alpha}$ and apply Theorem \ref{pom2} and Theorem \ref{Beckerthm2}.
\end{proof}

In \cite{Krzyz:1987}, a quasiconformal extension criterion for the class of strongly spirallike functions has been studied.
It shows that if there exists a $\beta \in [0,1)$ such that $e^{i\alpha}zf'(z)/f(z)$ lies in the sectoral domain $\{z : |\arg z| < \pi\beta/2\}$ for all $z \in \D$, then $f$ can be extended to a $k$-quasiconformal automorphism of $\CC$.
The dilatation of the extended mapping is implicitly given by $k = \sin (\pi\beta/2)$.
The same result was proved also in \cite[Korollar 5.3]{Gall:1986} though it does not mention the dilatation $k$.
A special case of Theorem \ref{spiralthm} when $\alpha = 0$ was provided in \cite{Betker:1992} as a consequence of an extension of Theorem \ref{Beckerthm}.

\

\subsection{Bazilevi\v c functions 1}

Let $\alpha>0$ and $\beta \in \R$.
In the following two subsections we deal with two quasiconformal extension criteria for a Bazilevi\v c function of type $(\alpha,\beta)$. 
Here, A function $f \in \A$ is called \textit{Bazilevi\v c of type $(\alpha, \beta)$} if
\begin{equation}\label{bazilevic}
f(z)=
\left[
(\alpha +i\beta)
\int_{0}^{z}
 g(\zeta)^{\alpha} h(\zeta) \zeta^{i\beta-1}d\zeta
\right]^{1/(\alpha+i \beta)}
\end{equation}
for a starlike univalent function $g \in \A$ and an analytic function $h$ with $h(0)=1$ satisfying $\textit{Re} (e^{i\lambda}h) > 0$ in $\D$ for some $\lambda \in \R$.
We denote by $\B(\alpha,\beta)$ the class of such functions.

If the above $h$ in \eqref{bazilevic} satisfies $|\arg e^{i\lambda} h| \leq \pi \gamma /2$ in $\D$ for some $\lambda \in [0,1] \hspace{5pt}(0 \leq \gamma \leq 1$ with $h(0)=1$, then we say that $f$ is a Bazilevi\v c function of order $\gamma$, and denote by $f \in \B(\alpha,\beta,\gamma)$.
This notion appears in \cite{Gall:1986}\footnote{The author would like to thank Professor Wolfram Koepf for his help under which the author could obtain a copy of the dissertation.} which is devoted to
the study of the class $\B(\alpha,\beta,\gamma)$ and contains the following quasiconformal extension criterion;
\textit{
Let $\alpha>0,\beta \in \R$ and $\gamma \in [0,1)$.
Then for $f \in \B(\alpha,\beta,\gamma)$ the following two conditions imply that $f(\round \D)$ is a rectifiable quasicircle, in particular, $f$ has a quasiconformal extension to $\CC$;
\def\labelenumi{\arabic{enumi}.}
\begin{enumerate}
\item $\dstyle \alpha < \frac{1-\gamma}{4}$,\\[-15pt]
\item There exists a starlike function $g \in \A$ in \eqref{bazilevic} such that
\begin{equation}\label{valueM}
\lim_{r \to 1} \left[ \frac{\log \max_{|z| = r}|g(z)|}{\log \frac{1}{1-r}}\right]
<
\frac{1-\gamma}{2\alpha}.
\end{equation}
\end{enumerate}
}

\no
It is known that the limit value of the left-hand side in \eqref{valueM} exists for all starlike functions in $\A$ (see e.g. \cite{pom:1962a}).
The above theorem does not estimate the dilatation of the extended quasiconformal mapping explicitly.

It follows from a result of Sheil-Small \cite{Sheil-Small:1972} that
$$
\Re 
  \left\{
    1 + \frac{zf''(z)}{f'(z)} + (\alpha + i \beta - 1)\frac{zf'(z)}{f(z)}
  \right\}
>0
$$
for all $z \in \D$ implies that $f \in \B(\alpha,\beta)$.
Now we refine this result to $k$-quasiconformal extension criterion (compare with \cite[Theorem 2]{Hotta:b});

\begin{thm}
Let $\alpha >0,\, \beta \in \R,\,k \in [0,1)$ and $f \in \A$.
If $f$ satisfies
\begin{equation}\label{convexstarlike}
  \left[
    1 + \frac{zf''(z)}{f'(z)} + (\alpha + i \beta - 1) \frac{zf'(z)}{f(z)}
  \right]
  \in U(k)
\end{equation}
for all $z \in \D$, then $f$ can be extended to a $k$-quasiconformal automorphism of $\C$.
\end{thm}

\begin{proof}
Let
$$
f_{t}(z) =  f(z)
  \left\{
  1+ (e^{t}-1)\frac{zf'(z)}{f(z)}
  \right\}^{1/(\alpha + i \beta)}.
$$
In that case $a_{1}(t) = f_{t}'(0) = e^{t/(\alpha + i \beta)}$ and therefore $\lim_{t \to \infty}|a_{1}(t)| = \infty$ since $\alpha > 0$. 
A straightforward calculation shows that
$$
\frac{1}{p(z,t)} =
\frac{1}{e^{t}}(\alpha + i \beta) +
\left(
  1 - \frac{1}{e^{t}}
\right) 
\left(
  1 + \frac{zf''(z)}{f'(z)} + (\alpha + i \beta - 1) \frac{zf'(z)}{f(z)}
\right).
$$
Since the assumption \eqref{convexstarlike} implies $\alpha + i \beta \in U(k)$ (consider the case $z = 0$), we have $1/p(z,t) \in U(k)$ by \eqref{convexstarlike}.
Consequently, it follows from Theorem \ref{pom2} and Theorem \ref{Beckerthm2} that $f$ has a $k$-quasiconformal extension to $\C$.
\end{proof}

\

\subsection{Bazilevi\v c functions 2}

A L\"owner chain for a Bazilevi\v c function is known \cite{Pom:1965} as
\begin{equation}\label{bazileL}
f_{t}(z) = 
\left[
(\alpha +i\beta)
\int_{0}^{z}
 g(\zeta)^{\alpha} \Big\{h(\zeta)+t h_{0}(\zeta)\Big\} \zeta^{i\beta-1}d\zeta
\right]^{1/(\alpha+i \beta)}
\end{equation}
where
$$
h_{0}(z) = i\beta + \alpha \frac{zg'(z)}{g(z)}.
$$
But we cannot apply Theorem \ref{Beckerthm} to \eqref{bazileL} because
$$
\frac{1}{p(z,t)}
= 
h(z) + t h_{0}(z)
$$
and hence it does not satisfy the assumption \eqref{becker} when $t$ tends to $\infty$.
We can avoid this obstacle by reparametrizing $t$ with $e^{t}-1$;

\begin{thm}
Let $\alpha>0,\,\beta \in \R$ and $k \in [0,1)$.
For $f \in \B (\alpha,\beta)$, we suppose that functions $h$ and $g$ of \eqref{bazilevic} satisfy $h(z) \in U(k)$ and $i\beta + \alpha (zg'(z)/g(z)) \in U(k)$ for all $z \in \D$, respectively.
Then $f$ can be extended to a $k$-quasiconformal automorphism of $\C$.
\end{thm}

\begin{proof}
Let
$$
f(z,t) = 
\left[
(\alpha +i\beta)
\int_{0}^{z}
 g(\zeta)^{\alpha} [h(\zeta)+(e^{t}-1)h_{0}(\zeta)] \zeta^{i\beta-1}d\zeta
\right]^{1/(\alpha+i \beta)}
$$
where
$$
h_{0}(z) = i\beta + \alpha \frac{zg'(z)}{g(z)}.
$$
Since
$$
\dot{f_{t}}(z) 
=
f(z,t)^{1-(\alpha+i\beta)}\cdot e^{t}g(z)^{\alpha}z^{i\beta}
$$
and
$$
zf_{t}'(z)
=
f(z,t)^{1-(\alpha+i\beta)}\cdot 
g(z)^{\alpha} [h(z)+(e^{t}-1)h_{0}(z)] z^{i\beta},
$$
we have
$$
\frac{1}{p(z,t)} = \frac{1}{e^{t}}h(z) + \left(1-\frac{1}{e^{t}}\right)h_{0}(z).
$$
We also have
$$
f_{t}'(0) = \left(1 + (e^{t} -1)(\alpha + i \beta)\right)^{1/(\alpha + i \beta)} 
$$
and therefore $\lim_{t \to \infty}|f_{t}'(0)| = \infty$.
Consequently, our assertion follows from Theorem \ref{pom2} and Theorem \ref{Beckerthm2}.
\end{proof}

\

\section*{Acknowledgement} 
The author is deeply grateful to Professor Toshiyuki Sugawa for many valuable suggestions and comments that have improved this paper. 
The author also thanks the referee for reading the manuscript carefully and giving several helpful comments.

\


\bibliographystyle{amsplain}
\bibliography{/bibdata.bib}
\end{document}